\documentclass[DIV=11,11pt]{scrartcl}
\usepackage[english]{babel}
\usepackage[utf8]{inputenc}
\usepackage[T1]{fontenc}

\author{Friedrich Martin Schneider, Daniel Borchmann}
\title{Topological Entropy of Formal Languages}
\date{\today}

\usepackage[scale=0.95]{tgpagella}
\usepackage[scale=0.92]{tgheros}
\usepackage[scaled=0.83]{beramono}
\usepackage{mathpazo}

\usepackage{mathtools,amssymb}
\usepackage[thmmarks,amsmath,hyperref]{ntheorem}
\usepackage{array}
\usepackage[colorlinks,citecolor=blue,hidelinks]{hyperref}
\usepackage{enumitem}
\setlist{noitemsep}
\setlist[enumerate]{label={\upshape\arabic*)}}
\usepackage{faktor}

\usepackage[backend=bibtex,style=numeric-comp,url=true]{biblatex}
\usepackage{cleveref}

\newtheorem{Theorem}{Theorem}[section]

\newtheorem{Corollary}[Theorem]{Corollary}
\newtheorem{Lemma}[Theorem]{Lemma}
\newtheorem{Proposition}[Theorem]{Proposition}

\theoremstyle{definition}
\theoremsymbol{\ensuremath{\triangle}}

\newtheorem{Definition}[Theorem]{Definition}

\theoremstyle{plain}
\theoremsymbol{\scalebox{1.5}[1.5]{$\diamond$}}
\theorembodyfont{\normalfont\upshape}
\newtheorem{Remark}[Theorem]{Remark}
\newtheorem{Example}[Theorem]{Example}

\theoremstyle{nonumberplain}
\theoremheaderfont{\normalfont\itshape}
\theorembodyfont{\normalfont\upshape}
\theoremsymbol{\ensuremath{\square}}
\newtheorem{Proof}{Proof}

\renewcommand{\epsilon}{\varepsilon}
\renewcommand{\phi}{\varphi}

\DeclareSymbolFont{bbold}{U}{bbold}{m}{n}
\DeclareSymbolFontAlphabet{\mathbbold}{bbold}

\DeclareMathOperator{\egr}{egr}
\DeclareMathOperator{\ind}{ind}

\newcommand{\set}[1]{\{\,#1\,\}}
\newcommand{\abs}[1]{\lvert#1\rvert}
\newcommand{\subsets}[1]{\mathop{\mathfrak{P}}(#1)}
\newcommand{\NN}{\mathbb{N}}
\newcommand*{\im}[1]{\operatorname{im}(#1)}

\begin{document}

\maketitle

\begin{abstract}
  We introduce the notion of \emph{topological entropy} of a formal languages as
  the topological entropy of the minimal topological automaton accepting it.
  Using a characterization of this notion in terms of approximations of the
  Myhill-Nerode congruence relation, we are able to compute the topological
  entropies of certain example languages.  Those examples suggest that the
  notion of a \enquote{simple} formal language coincides with the language
  having zero entropy.
\end{abstract}

\section{Introduction}
\label{sec:introduction}

The Chomsky hierarchy classifies formal languages in levels of growing
complexity.  At its bottom it puts the class of regular languages, followed by
context-free and context-sensitive languages.  At the top of the hierarchy it
lists the class of all decidable languages.  As such, the Chomsky hierarchy
gives a method to assign a measure of \emph{complexity} to formal languages.

However, using the Chomsky hierarchy as a mean to asses the complexity of a
language has certain drawbacks.  The most severe drawback is that the
classification of the Chomsky hierarchy depends on a particular choice of
computation models, namely finite automata, non-deterministic pushdown-automata,
linear bounded automata, and Turing machines, respectively.  It can be argued
that this choice results in some contra-intuitive classifications: of course,
accepting a language to be \enquote{simple} as soon as it is accepted by a
finite automaton is reasonable.  The converse however is not: not every language
that cannot be accepted by a finite automaton is necessarily
\enquote{complicated}.

An example is the Dyck language $D$ with one sort of
parentheses~\cite{book/HoFL/ContextFreeLanguages}.  This is the language of all
words of balanced parentheses like (()()) and ((())) but not (())) or (.  This
language is context-free but not regular, and thus a classification by the
Chomsky hierarchy would make this language $D$ appear to be not so
\enquote{simple}.  On the other hand, there is a very simple machine model
accepting $D$, namely a two-state automaton with only one counter.  It is
reasonable to say that this kind of automaton is intuitively simple.  The
Chomsky hierarchy does not capture this: it puts the Dyck language with one sort
of parentheses in the same class as much more complicated languages like
palindromes.  And there even exist context-sensitive languages that can be
accepted by finite automata with only one counter.

To assess the complexity of a language one could now proceed as follows: given a
language $L$, what is the simplest form of computation model that is required to
accept $L$?  It is clear that this approach heavily depends on the notion of
\enquote{simplest computation model} and the fact that there is such one.
Indeed, it requires a hierarchy of all conceivable computation models to make
this approach work, an assumption that is hardly realizable.

Instead of considering all possible computation models, we propose another
approach, namely to consider \emph{one} computation model that works \emph{for
  every} language.  Then given a formal language $L$ one could ask what the
\enquote{simplest} instances of this particular computation model is that is
required to accept the given language $L$.  This then can be used to assign to
$L$ a measure of complexity that does not depend on a particular a-priori choice
of certain computational models.

More precisely, we shall show in this work that we can use the notion of
\emph{topological automata}~\cite{journals/corr/Steinberg13} to assign to every
formal languages a notion of \emph{entropy} that naturally reflects the
complexity of the formal languages.  As such, we make use of the following
facts: for every formal language there exists a topological automaton accepting
it.  Furthermore, for each topological automaton there exists a natural notion
of a \emph{smallest} automaton accepting the same language.  Finally, as
topological automata are a particular form of \emph{dynamical systems}, we can
naturally assign a measure of complexity to every topological automaton, namely
its \emph{entropy}.  Therefore, we can define the complexity of a formal
language $L$ as the \emph{entropy of the minimal topological automaton accepting
  it}.  We call this notion the \emph{topological entropy} of $L$.  Intuitively,
the lower the topological entropy of $L$ the simpler it is.  Languages with
vanishing entropy are thus the simplest of all formal languages.

An advantage of this approach is that it works for every formal language, and is
thus independent of a particular choice of computation models.  On the other
hand, one could argue that this approach is purely theoretical, as it may not
allow us to compute the entropy of formal languages easily.  However, we shall
show that it is indeed possible to compute the topological entropy for certain
examples of languages.  For this we use a characterization of the topological
entropy in terms of approximations of its Myhill-Nerode congruence relation.
Using this, it is not hard to show that all regular languages have entropy 0.
Moreover, we shall show that the Dyck language with one sort of parenthesis has
also entropy 0.  Both of these can thus be called \enquote{simple}, and
intuitively they are.  On the other hand, we shall also show that languages like
palindromes or Dyck languages with multiple sorts of parentheses do not have
zero entropy.

The paper is structured as follows.  We first introduce the notions of
\emph{topological automata} and \emph{entropy of semigroup actions} and formally
define the notion of \emph{topological entropy} of formal languages.  The main
part of this paper is then devoted to proof a characterization of topological
entropy that allows for a comparably easy way to compute it.  This is done in
\Cref{sec:char-topol-entr}.  We compute the entropy of some example languages in
\Cref{sec:examples}.  We also provide a characterization of the topological
entropy in terms of the \emph{entropic dimension} of suitable pseudo-ultrametric
spaces.  Finally, we shall summarize the results of this paper and sketch an
outline of future work.

\section{Topological Entropy of Formal Languages}
\label{sec:topol-autom-entr}

A variety of notions has been developed to assess different aspects of
complexity of formal languages.  Most of these notions have been devised with an
understanding of complexity in mind that comes with classical complexity theory,
and thus these notions are formulated as decision problems.  Examples for this
are the word problem and the equivalence problem for formal languages, and the
complexity of the formal languages is measured by the complexity class for which
these problems are complete.  Other notions quantify complexity by other means.
Examples are the \emph{state complexity}~\cite{Yu00statecomplexity} of a regular
language, which gives the complexity of the language as the number of states in
its minimal automaton, or the \emph{syntactic complexity} of a regular language,
which instead considers the size of corresponding syntactic
semigroup~\cite{conf/dlt/BrzozowskiY11}.

The core idea of the present article is to expand the methods of measuring a
formal language's complexity by a topological approach in terms of
\emph{topological entropy}, which proved tremendously useful to dynamical
systems. Topological entropy was introduced by Adler et
al. \cite{journal/tams/AdlerKMcA65} for single homeomorphisms (or continuous
transformations) on a compact Hausdorff space. The literature provides several
essentially different extensions of this concept for continuous group and
semigroup actions. Among others, there is an approach towards topological
entropy for continuous actions of finitely generated (pseudo-)groups due to Ghys
et al.~\cite{Ghys} (see also~\cite{Bis2,BisWalczak,Bis3,Walczak}), which has
also been investigated for continuous semigroup actions
in~\cite{Bis1,BisUrbanski,specification}.

By a dynamical system we mean a continuous semigroup action on a compact
Hausdorff topological space. Topological entropy measures the ability of an
observer to distinguish between points of the dynamical system just by
recognizing transitions \emph{at equal time intervals}, i.e.,~with respect to a
fixed generating system of transformations, starting from the initial
state. Since the above notion of dynamical system may very well be regarded as
the topological counterpart of a finite automaton, it seems natural to utilize
the dynamical approach for applications to automata theory.

To link dynamical systems to formal languages we shall use the already mentioned
notion of a \emph{topological automaton}~\cite{journals/corr/Steinberg13}.  This
notion has been introduced as a topological generalization of the usual notion
of a \emph{finite automaton} by allowing to have an infinite state space.
Indeed, topological automata share certain properties with finite automata.  For
example, for each topological automaton there exists a \emph{minimal}
topological automaton accepting the same language.  However, in contrast to
finite automata, \emph{every} formal language is accepted by some
\emph{topological automaton}, and not only regular languages.  This allows us to
uniformly treat all formal languages with one computation model.

Recall that a \emph{(deterministic) automaton} over an alphabet $\Sigma$ is a
tuple $\mathcal{A} = (Q, \Sigma,\delta,q_{0},F)$ consisting of a finite set $Q$
of \emph{states}, a \emph{transition function} $\delta \colon Q \times \Sigma
\to Q$, a set $F \subseteq Q$ of \emph{final} states, and an \emph{initial
  state} $q_{0} \in Q$.  The transition function is usually extended to the set
of all words over $\Sigma$ by virtue of
\begin{align*}
  \delta^{*}(q, \epsilon) &\coloneqq q, \\
  \delta^{*}(q, wa) &\coloneqq \delta(\delta(q, w), a)
\end{align*}
for $q \in Q$, $a \in \Sigma$, and $w \in \Sigma^{*}$. The
\emph{language accepted} by $\mathcal{A}$ is then
\begin{equation*}
  \mathcal{L}(\mathcal{A}) \coloneqq \{ w \in \Sigma^{\ast} \mid \delta^{\ast}
  (q_{0},w) \in F \}.
\end{equation*}
It is not hard to see that the function $\delta^{*}$ is a monoid action of
$\Sigma^{*}$ on $Q$ -- indeed it is the unique monoid action of $\Sigma^{*}$ on
$Q$ extending $\delta$.

The notion of deterministic finite automata can now be extended to an infinite
state set as follows.  Throughout this article a \emph{continuous action} of a
semigroup or monoid $S$ on a topological space $X$ is an action $\alpha$ of $S$
on $X$ such that $\alpha_{s}\colon X \to X, \, \alpha_{s}(x) = \alpha(x, s)$ is
continuous for every $s \in S$.  Note that the latter just means that $\alpha
\colon X \times S \to X$ is continuous where $S$ is endowed the discrete
topology.

\begin{Definition}
  A \emph{topological automaton over an alphabet $\Sigma$} is a tuple
  $\mathcal{A} = (X,\Sigma,\alpha, x_{0}, F)$ consisting of
  \begin{itemize}
  \item a compact Hausdorff space $X$, called the set of \emph{states} of
    $\mathcal{A}$
  \item a continuous action $\alpha$ of $\Sigma^{\ast}$ on $X$, called the
    \emph{transition function} of $\mathcal{A}$,
  \item a point $x_{0} \in X$, called the \emph{initial state} of $\mathcal{A}$,
    and
  \item a clopen subset $F \subseteq X$, called the set of \emph{final states}
    of $\mathcal{A}$.
  \end{itemize}
  We say that $\mathcal{A}$ is \emph{trim} if $\alpha (x_{0},\Sigma^{\ast})$ is
  dense in $X$. The \emph{language recognized by $\mathcal{A}$} is defined as
  \begin{equation*}
    \mathcal{L}(\mathcal{A}) \coloneqq \{ w \in \Sigma^{\ast} \mid \alpha (x_{0}, w ) \in F \} .
  \end{equation*}

  Let $\mathcal{B} = (Y, \Sigma, \beta, y_{0}, G)$ be another topological
  automaton.  We shall say that $\mathcal{A}$ and $\mathcal{B}$ are
  \emph{isomorphic}, and write $\mathcal{A} \cong \mathcal{B}$, if there exists
  a homeomorphism $\phi \colon X \to Y$ such that
  \begin{equation*}
    \phi(\alpha(x,\sigma)) = \beta(\phi(x), \sigma)
  \end{equation*}
  for all $x \in X, \sigma \in \Sigma$, $\phi(x_{0}) = y_{0}$, and $\phi(F) = G$.
\end{Definition}

Evidently, isomorphic automata accept the same language.

Observe that every automaton accepting $L$ can be turned into an automaton that
is trim: if $\mathcal{A} = (X, \Sigma, \alpha, x_{0}, F)$ is a topological
automaton accepting $L$, then replacing $X$ with $\overline{\alpha(x_{0},
  \Sigma^{*})}$ and $F$ with $F \cap \overline{\alpha(x_{0}, \Sigma^{*})}$
always yields a trim automaton accepting the same language $L$.

As already stated, and in contrast to regular languages, every formal language
$L \subseteq \Sigma^{*}$ is accepted by a topological automaton,
cf.~\cite[Proposition~2.1]{journals/corr/Steinberg13}.

\begin{Proposition}
  \label{prop:topological-automata-accepting-L-always-exist}
  Let $L \subseteq \Sigma^{*}$ and $\chi_{L}$ the characteristic function of
  $L$.  Equip $X \coloneqq \set{ 0, 1 }^{\Sigma^{*}}$ with the product topology,
  and define the mapping $\delta \colon X \times \Sigma^{*} \to X$ by
  \begin{equation*}
    \delta(f, u)(v) \coloneqq f(uv).
  \end{equation*}
  Then $L$ is accepted by the topological automaton $(X, \Sigma, \delta,
  \chi_{L}, T)$ for $T \coloneqq \set{ f \in X \mid f(\epsilon) = 1 }$.
\end{Proposition}

With the notation of \Cref{prop:topological-automata-accepting-L-always-exist},
we define the \emph{minimal automaton of $L$} to be
\begin{equation*}
  \mathcal{A}_{L} = \bigl(\overline{\chi_{L}(\Sigma^{*})}, \Sigma, \delta, \chi_{L}, T_{L} \bigr),
\end{equation*}
where $\overline{\chi_{L}(\Sigma^{*})}$ is the closure of $\chi_{L}(\Sigma^{*})$
in $\set{0, 1}^{\Sigma^{*}}$, and $T_{L} = T \cap
\overline{\chi_{L}(\Sigma^{*})}$.  Clearly, $\mathcal{A}_{L}$ is trim.  Indeed
we have the following fact that justifies to call $\mathcal{A}_{L}$ minimal,
cf.~\cite[Theorem~2.2]{journals/corr/Steinberg13}.

\begin{Proposition}
  Let $L \subseteq \Sigma^{*}$, and let $\mathcal{A} = (X, \Sigma, x_{0},
  \delta, F)$ be a topological automaton accepting $L$.  Then $\mathcal{A} \cong
  \mathcal{A}_{L}$ if and only if for every trim automaton $\mathcal{B} = (Y,
  \Sigma, y_{0}, \lambda, G)$ accepting $L$ there exists a uniquely determined
  surjective continuous function $\phi \colon Y \to X$ satisfying
  $\phi(\lambda(y, \sigma)) = \delta(\phi(y), \sigma)$ and $\phi(y_{0}) =
  x_{0}$.  Moreover, in this case the unique $\phi$ satisfies $G =
  \phi^{-1}(F)$.
\end{Proposition}

Since $\mathcal{A}_{L} \cong \mathcal{A}_{L}$, this proposition immediately
yields that the minimal automaton is indeed minimal in the above sense.
Moreover, in the case that $L$ is regular, $\mathcal{A}_{L}$ is finite and is
the usual minimal automaton of regular languages.

\begin{Example}
  Let $\Sigma$ be a finite alphabet and let $a,b \in \Sigma$, $a \ne b$.  We
  consider the \emph{Alexandroff compactification} $\mathbb{Z}_{\infty}$ of the
  discrete space of integers $\mathbb{Z}$, that is the set $\mathbb{Z}_{\infty}
  = \mathbb{Z} \cup \{ \infty \}$ equipped with the topology
  \begin{equation*}
    \{ M \subseteq \mathbb{Z} \cup \{ \infty\} \mid \infty \in M \Longrightarrow
    \mathbb{Z} \setminus M \text{ is finite} \} .
  \end{equation*}
  We define an action $\alpha$ of $\Sigma^{\ast}$ on $\mathbb{Z}_{\infty}$ by
  setting $\alpha (m,a) = m+1$, $\alpha (m,b) = m-1$ and $\alpha (m,c) = m$ for
  all $m \in \mathbb{Z}_{\infty}$ and $c \in \Sigma\setminus \{ a,b \}$.  Then
  $\alpha$ constitutes a continuous action of $\Sigma^{\ast}$ on
  $\mathbb{Z}_{\infty}$, and for each $n \in \mathbb{N}$ the topological
  automaton $\mathcal{A} = (\mathbb{Z}_{\infty},\Sigma, \alpha,0,\{ n \})$
  accepts the language $L = \bigl\{ w \in \Sigma^{\ast} \bigm| \vert w \vert_{a}
  = \vert w \vert_{b} + n \bigr\}$.
\end{Example}

We now shall express the complexity of the language $L$ accepted by a
topological automaton $\mathcal{A} = (Q, \Sigma, \alpha, x_{0}, F)$ by the
\emph{topological entropy} of the continuous action $\alpha$ of $\Sigma^{*}$ on
$Q$~\cite{journals/koreamath/XinhuaL08,journal/tams/AdlerKMcA65,journal/etds/BlanchardHM00}.
To this end, we shall first fix some useful notation and recall some important
definitions about continuous actions on compact Hausdorff spaces.

Let $X$ again be a compact Hausdorff space. We shall denote by $\mathcal{C}(X)$
the set of all finite open covers of $X$. If $f \colon X \to X$ is continuous
and $\mathcal{U} \in \mathcal{C}(X)$, then $f^{-1}(\mathcal{U}) := \{ f^{-1}(U)
\mid U \in \mathcal{U} \}$ is a finite open cover of $X$ as well. Given
$\mathcal{U}, \mathcal{V} \in \mathcal{C}(X)$, we say that $\mathcal{V}$
\emph{refines} $\mathcal{U}$ and write $\mathcal{U} \preceq \mathcal{V}$ if
\begin{equation*}
  \forall V \in \mathcal{V} \, \exists U \in \mathcal{U} \colon \, V \subseteq U ,
\end{equation*}
and we say that $\mathcal{U}$ and $\mathcal{V}$ are \emph{refinement-equivalent} and write
$\mathcal{U} \equiv \mathcal{V}$ if $\mathcal{U} \preceq \mathcal{V}$ and
$\mathcal{V} \preceq \mathcal{U}$. Furthermore, if $(\mathcal{U}_{i} \mid i \in I)$ is a
finite family of finite open covers of $X$, then
\begin{equation*}
  \bigvee_{i\in I} \mathcal{U}_{i} \coloneqq \bigl\{\, \bigcap_{i\in I} U_{i} \bigm|
  (U_{i})_{i\in I} \in \prod_{i\in I} \mathcal{U}_{i} \,\bigr\}.
\end{equation*}
is a finite open cover of $X$ as well.  For $\mathcal{U} \in \mathcal{C}(X)$ let
\begin{equation*}
  N(\mathcal{U}) \coloneqq \inf\bigl\{\, \abs{\mathcal{V}} \bigm| \mathcal{V} \subseteq
  \mathcal{U}, X = \bigcup \mathcal{V} \,\bigr\}.
\end{equation*}

In preparation for some later considerations, let us recall the following basic
observations.

\begin{Remark}[\cite{journal/tams/AdlerKMcA65}]
  \label{rem:entropy}
  Let $X$ be a compact Hausdorff space, $\mathcal{U}, \mathcal{V} \in
  \mathcal{C}(X)$, $I$ be a finite set, $(\mathcal{U}_{i})_{i \in I},
  (\mathcal{V}_{i})_{i \in I} \in \mathcal{C}(X)^{I}$, and $f \colon X \to X$ be
  a continuous map. Then the following statements hold:
  \begin{enumerate}
  \item $\mathcal{U} \preceq \mathcal{V} \, \Longrightarrow \, N(\mathcal{U})
    \leq N(\mathcal{V})$,
  \item $\mathcal{U} \preceq \mathcal{V} \, \Longrightarrow \,
    f^{-1}(\mathcal{U}) \preceq f^{-1}(\mathcal{V})$,
  \item $(\forall i \in I \colon \, \mathcal{U}_{i} \preceq \mathcal{V}_{i}) \,
    \Longrightarrow \, \bigvee_{i \in I} \mathcal{U}_{i} \preceq \bigvee_{i \in
      I} \mathcal{V}_{i}$.
  \end{enumerate}
\end{Remark}

Now we come to dynamical systems, i.e., continuous semigroup actions. Let $S$ be
a semigroup and consider a continuous action $\alpha$ of $S$ on $X$.  For
$\mathcal{U} \in \mathcal{C}(X)$ we write
\begin{equation*}
  s^{-1}(\mathcal{U}) \coloneqq \alpha_{s}^{-1}(\mathcal{U}).
\end{equation*}
For every finite $F \subseteq S$ and $\mathcal{U} \in \mathcal{C}(X)$ let
\begin{equation*}
  (F : \mathcal{U})_{\alpha} \coloneqq N\bigl(\bigvee_{s \in F} s^{-1}(\mathcal{U})\bigr).
\end{equation*}
Assume $F$ to be a finite generating subset of $S$. If $\mathcal{U}$ is a finite
open cover of $X$, then we define
\begin{equation*}
  \eta (\alpha ,F, \mathcal{U}) := \limsup_{n \to \infty} \frac{\log_{2}
    (F^{n}:\mathcal{U})_{\alpha}}{n} .
\end{equation*}
Furthermore, the \emph{topological entropy of $\alpha$ with respect to $F$} is
defined to be the quantity
\begin{equation*}
  \eta (\alpha ,F) := \sup \{ \eta(\alpha, F,\mathcal{U}) \mid \mathcal{U} \in \mathcal{C}(X) \} .
\end{equation*}
Of course, the precise value of this quantity depends on the choice of a finite
generating system. However, we observe the following fact.

\begin{Proposition}
  Let $S$ be a semigroup and let $\alpha$ be a continuous action of $S$ on some
  compact Hausdorff space $X$. Suppose $E,F \subseteq S$ to be finite subsets
  generating $S$. Then
  \begin{equation*}
    \frac{1}{m} \cdot \eta(\alpha,F) \leq \eta(\alpha,E) \leq n \cdot \eta(\alpha,F) ,
  \end{equation*}
  where $m := \inf \{ k \in \mathbb{N} \mid F \subseteq E^{k} \}$ and
  $n := \inf \{ k \in \mathbb{N} \mid E \subseteq F^{k} \}$.
\end{Proposition}

\begin{Proof}
  Let $\mathcal{U} \in \mathcal{C}(X)$. Evidently, $(E^{k}:\mathcal{U}) \leq
  (F^{kn}:\mathcal{U})$ for all $k \in \mathbb{N}$, whence
  \begin{align*}
    \eta(\alpha ,E,\mathcal{U})
    &= \limsup_{k \to \infty} \frac{\log_{2} (E^{k}:\mathcal{U})_{\alpha}}{k}
      \leq \limsup_{k \to \infty} \frac{\log_{2} (F^{kn}:\mathcal{U})_{\alpha}}{k} \\
    &= n\limsup_{k \to \infty} \frac{\log_{2} (F^{kn}:\mathcal{U})_{\alpha}}{kn} \leq n\limsup_{k
      \to \infty} \frac{\log_{2} (F^{k}:\mathcal{U})_{\alpha}}{k}
      = n \cdot \eta(\alpha ,F,\mathcal{U})
  \end{align*}
  Thus, $\eta(\alpha ,E,\mathcal{U}) \leq n \cdot \eta(\alpha ,F,\mathcal{U})$. This shows
  that $\eta(\alpha, E) \leq n\eta(\alpha, F)$. Due to symmetry, it follows that
  $\eta(\alpha, F) \leq m \cdot \eta(\alpha, E)$ as well.
\end{Proof}

With all the necessary notions in place we are finally able to define our notion
of entropy of formal languages.

\begin{Definition}
  Let $L \subseteq \Sigma^{*}$ and let $\mathcal{A}_{L} = (X, \Sigma, x_{0},
  \delta, F)$ be the minimal automaton of $L$.  Then the \emph{entropy} of $L$
  is the entropy $\eta(\delta, \Sigma \cup \set{ \epsilon })$ of $\delta$ with
  respect to $\Sigma \cup \set{ \epsilon }$.
\end{Definition}

\section{A Characterization}
\label{sec:char-topol-entr}

We claim that the definition of topological entropy is natural.  Yet
\emph{computing} using the definition alone may not work very well.  It is the
purpose of this section to remedy this issue by providing an alternative
characterization of the topological entropy of formal languages.  For this we
exploit another way of considering formal languages as dynamical systems.

To view a formal language $L$ over an alphabet $\Sigma$ as some kind of
dynamical system we take inspiration from the characterization of regular
languages as languages whose Myhill-Nerode congruence relation $\Theta(L)$ has
finite index.  Recall that for $u, v \in \Sigma$ we have
\begin{equation*}
  (u,v) \in \Theta(L) \iff \forall w \in \Sigma^{*}
  \colon (uw \in L \iff vw \in L).
\end{equation*}
The relation $\Theta(L)$ can be seen as some way of measuring the complexity of
$L$: if $L$ is regular, the number of equivalence classes is finite and equals
the number of states in the minimal automaton of $L$. Indeed, this is the idea
behind the notion of state complexity.

However, if $L$ is not regular this measure is not available anymore.  We shall
remedy this by not considering the \emph{number} of equivalence classes of
$\Theta(L)$, but by considering the \emph{growth} of the number of equivalence
classes of a particular approximation of $\Theta(L)$.  Based on this growth we
introduce our characterization of topological entropy of $L$.

Let us first recall some basic notation.  Let $\Theta$ be an equivalence
relation on a set $Y$.  For $y \in Y$ we put $[y]_{\Theta} := \set{ x \in Y \mid
  (x, y) \in \Theta}$.  Then $Y / \Theta := \set{ [y]_{\Theta} \mid y \in Y}$.
Furthermore, the \emph{index} of $\Theta$ on $Y$ is defined as $\ind(\Theta) :=
\abs{Y / \Theta}$.  For a mapping $f \colon X \to Y$ we set $f^{-1}(\Theta)
\coloneqq \set{(s,t) \in X\times X \mid (f(x),f(y)) \in \Theta}$.  Clearly,
$f^{-1}(\Theta)$ then constitutes an equivalence relation on $X$.

Now let $\Sigma$ be an alphabet.  The \emph{Nerode congruence} of a language $L
\subseteq \Sigma^{\ast}$ is the equivalence relation
\begin{equation*}
  \Theta (L) \coloneqq \{ (u,v) \in \Sigma^{\ast} \times \Sigma^{\ast} \mid \forall w
  \in \Sigma^{\ast} \colon uw \in L \Leftrightarrow vw \in L \} .
\end{equation*}
Recall that $L$ is regular if and only if it is accepted by an automaton.  The
following characterization of regular languages in terms of the Nerode
congruence relation is well-known.
\begin{Theorem}[Myhill-Nerode]
  \label{Theorem:Myhill-Nerode}%
  Let $\Sigma$ be a finite alphabet. A language $L \subseteq \Sigma^{\ast}$ is
  regular if and only if $\Theta (L)$ has finite index.
\end{Theorem}

Starting from this characterization we shall now make precise what we mean by
\emph{approximating} the relation $\Theta$.  For this we introduce another type
of equivalence relation.

\begin{Definition}
  Let $\Sigma$ be an alphabet.  For $F \subseteq \Sigma^{\ast}$ finite and $L
  \subseteq \Sigma^{\ast}$ define the function $\chi_{F, L} \colon \Sigma^{*}
  \to \set{0,1}^{F}$ by
  \begin{equation*}
    \chi_{F, L}(u)(w) \coloneqq
    \begin{cases}
      1 & \text{if } uw \in L \\
      0 & \text{otherwise}
    \end{cases}
  \end{equation*}
  for $u \in \Sigma^{*}$ and $w \in F$.  Now let
  \begin{equation*}
    \Theta(F, L) \coloneqq \ker(\chi_{F, L}).
  \end{equation*}

\end{Definition}

Now, the equivalence relations $\Theta(F, L)$ constitute an approximation of
$\Theta(L)$ in the sense that
\begin{equation}
  \label{eq:1}
  \Theta (L) = \bigcap \set{ \Theta(F,L) \mid F \subseteq \Sigma^{*} \text{ finite} }.
\end{equation}
Note that $\ind\Theta(F, L) = \lvert \im{\chi_{F, L}} \rvert \leq 2^{\lvert F
  \rvert}$.  In particular, $\Theta(F, L)$ has finite index, and thus the
following definition is reasonable.

Therefore, as $(\Theta(\Sigma^{(n)}, L) \mid n \in \mathbb N)$ may be regarded
as an approximation of the Myhill-Nerode congruence $\Theta(L)$, it seems
natural to consider the exponential growth rate of the corresponding index
sequence as a measure of complexity for a given formal language $L$.  This
motivates the following definition.

\begin{Definition}
  Let $\Sigma$ be an alphabet, and denote with $\mathcal{F}(\Sigma^{*})$ the set
  of finite subsets of $\Sigma^{*}$.  Define
  \begin{equation*}
    \gamma \colon \mathcal{F}(\Sigma^{\ast}) \times \mathcal{P}(\Sigma^{\ast})
    \to \mathbb{N}, \, (F,L) \mapsto \ind \Theta (F,L).
  \end{equation*}
  Given $L \subseteq \Sigma^{\ast}$, we call
  \begin{equation*}
    \gamma_{L} \colon \mathcal{F}(\Sigma^{\ast}) \to \mathbb{N},
    \, F \mapsto \gamma (F,L)
  \end{equation*}
  the \emph{Myhill-Nerode complexity function of $L$}. The \emph{Myhill-Nerode
    entropy} of a language $L \subseteq \Sigma^{\ast}$ is defined to be
  \begin{equation*}
    h(L) \coloneqq \limsup_{n \to \infty} \frac{\log_{2} \gamma_{L}(\Sigma^{(n)})}{n},
  \end{equation*}
  where $\Sigma^{(n)}$ is the set of all words over $\Sigma$ of length at most $n$.
\end{Definition}

The Myhill-Nerode complexity function of $L$ has some immediate properties that
we collect in the next proposition.

\begin{Proposition}
  Let $\Sigma$ be a finite alphabet, let $E,F \subseteq \Sigma^{\ast}$ be
  finite, and let $L,L_{0},L_{1} \subseteq \Sigma^{\ast}$.  Then
  \begin{enumerate}
  \item $\gamma(F, \emptyset) = \gamma(F, \Sigma^{*}) = 1$, and thus
    $h(\emptyset) = h(\Sigma^{*}) = 0$.
  \item $\gamma(E \cup F, L) \leq \gamma(E, L) \cdot \gamma(F, L)$.  If $E \subseteq F$, then
    $\gamma(E, L) \leq \gamma(F, L)$.
  \item $\gamma(F, L) = \gamma(F, \Sigma^{*} \setminus L)$, and hence
    $h(L) = h(\Sigma^{*} \setminus L)$.
  \item $\gamma(F, L_{0} \cup L_{1}) \leq \gamma(F, L_{0}) \cdot \gamma(F, L_{1})$, and thus
    $h(L_{0} \cup L_{1}) \leq h(L_{0}) + h(L_{1})$.
  \item $\gamma(F, L_{0} \cap L_{1}) \leq \gamma(F, L_{0}) \cdot \gamma(F, L_{1})$, and thus
    $h(L_{0} \cap L_{1}) \leq h(L_{0}) + h(L_{1})$.
  \end{enumerate}
\end{Proposition}

The goal of the rest of this section is now to show that the Myhill-Nerode
complexity of $L$ coincides with the topological entropy of $L$.  We start this
endeavor by showing that the Myhill-Nerode entropy of a formal language is
bounded from above by the entropy of any topological automaton accepting it. In
the case that the automaton is trim, these two notions even coincide.

\begin{Theorem}
  \label{thm:entropy-of-automaton-as-entropy-of-cover}
  Suppose $\mathcal{A} = (X,\Sigma,\alpha,x_{0},F)$ to be a topological
  automaton. Consider $S := \Sigma \cup \{ \epsilon \}$ and $\mathcal{U} := \{
  F, \, X\setminus F \}$. Then $h(L(\mathcal{A})) \leq \eta
  (\alpha,S,\mathcal{U})$. If $\mathcal{A}$ is trim, then $h(L(\mathcal{A})) =
  \eta (\alpha,S,\mathcal{U})$.
\end{Theorem}

We prove this theorem with the following three auxiliary statements.

\begin{Lemma}
  \label{lem:entropy-of-automaton-as-entropy-of-cover}
  Let $\mathcal{A} = (X,\Sigma,\alpha,x_{0},F)$ be a topological automaton. Let
  $\Phi \colon \Sigma^{\ast} \to X, \, w \mapsto \alpha(x_{0},w)$ and
  $\mathcal{U} := \{ F, \, X\setminus F \}$. Consider a finite subset
  $E \subseteq \Sigma^{\ast}$ as well as the equivalence relation
  \begin{equation*}
    \Lambda_{E} := \{ (x,y) \mid \forall w \in E \colon \, \alpha (x,w) \in F \Longleftrightarrow
    \alpha (y,w) \in F\} .
  \end{equation*}
  Then the following statements hold:
  \begin{enumerate}
  \item
    $X/\Lambda_{E} = \left( \bigvee\nolimits_{w \in E} w^{-1}(\mathcal{U}) \right) \setminus \{
    \emptyset \}$.
  \item $\Theta (E,L(\mathcal{A})) = (\Phi \times \Phi)^{-1}(\Lambda_{E})$.
  \item If $\mathcal{A}$ is trim, then $\Phi (\Sigma^{\ast}) \cap V \ne \emptyset$ for every
    $V \in X/\Lambda_{E}$.
  \end{enumerate}
\end{Lemma}

\begin{Proof}
  (1): We observe that
  $\mathcal{V} := \left( \bigvee\nolimits_{w \in E} w^{-1}(\mathcal{U}) \right) \setminus
  \{ \emptyset \}$
  constitutes a finite partition of $X$ into clopen subsets. For any $V \in \mathcal{V}$
  and $x \in V$, we observe that
  \begin{align*}
    [x]_{\Lambda_{E}}
    &= \{ y \in Y \mid \forall w \in E \colon \alpha (x,w) \in F \Leftrightarrow \alpha
      (y,w) \in F \} \\
    &= \{ y \in Y \mid \forall w \in E \colon x \in w^{-1}(F) \Leftrightarrow y \in w^{-1}(F) \} \\
    &= \{ y \in Y \mid \forall w \in E\, \forall U \in \mathcal{U}\, \colon x \in
      w^{-1}(U) \Leftrightarrow y \in w^{-1}(U) \} \\
    &= \{ y \in Y \mid \forall W \in \mathcal{V} \colon x \in W \Leftrightarrow y \in W \} \\
    &= \{ y \in Y \mid y \in V \} \\
    &= V
  \end{align*}
  We conclude that $X/\Lambda_{E} = \mathcal{V}$.

  (2): Let $L := L(\mathcal{A})$. For any two words $u,v \in \Sigma^{\ast}$, it follows that
  \begin{align*}
    (u,v) \in \Theta(E,L) \
    &\Longleftrightarrow \ \forall w \in E \colon uw \in L \Leftrightarrow vw \in L \\
    &\Longleftrightarrow \ \forall w \in E \colon \alpha (x_{0},uw) \in F \Leftrightarrow
      \alpha (x_{0},vw) \in F \\
    &\Longleftrightarrow \ \forall w \in E \colon \alpha (\alpha (x_{0},u),w) \in F
      \Leftrightarrow \alpha (\alpha (x_{0},v),w) \in F \\
    &\Longleftrightarrow \ \forall w \in E \colon \alpha (\Phi(u),w) \in F \Leftrightarrow
      \alpha (\Phi(v),w) \in F \\
    &\Longleftrightarrow \ (\Phi(u),\Phi(v)) \in \Lambda_{E} .
  \end{align*}
  That is, $\Theta (E,L) = (\Phi \times \Phi)^{-1}(\Lambda_{E})$.

  (3): By (1), the set $X/\Lambda_{E}$ is a collection of open, non-empty subsets of
  $X$. If $\mathcal{A}$ is trim, then $\Phi (\Sigma^{\ast})$ is dense in $X$, and thus
  $\Phi(\Sigma^{*}) \cap V \neq \emptyset$ for every $V \in X/\Lambda_{E}$.
\end{Proof}

\begin{Proposition}
  \label{prop:entropy-of-automaton-as-entropy-of-cover}
  Let $\mathcal{A} = (X,\Sigma,\alpha,x_{0},F)$ be a topological automaton and let
  $\mathcal{U} := \{ F, \, X\setminus F \}$. Consider a finite subset $E \subseteq
  \Sigma^{\ast}$.
  Then $\gamma_{L(\mathcal{A})}(E) \leq (E : \mathcal{U})_{\alpha}$. Furthermore, if $\mathcal{A}$
  is trim, then $\gamma_{L(\mathcal{A})}(E) = (E : \mathcal{U})_{\alpha}$.
\end{Proposition}

\begin{Proof}
  Let $L \coloneqq L(\mathcal{A})$ and $\mathcal{V} \coloneqq \bigvee\nolimits_{w \in
    E}(w^{-1}(\mathcal{U}))$. Since $\mathcal{V}\setminus \{ \emptyset \}$ constitutes a
  finite partition of $X$ into clopen subsets, $\mathcal{V} \setminus \set{\emptyset}$
  does not admit any proper subcover. Consequently, $N(\mathcal{V}) = \abs{\mathcal{V}
    \setminus \set{\emptyset}}$.  Applying
  \Cref{lem:entropy-of-automaton-as-entropy-of-cover}, we conclude
  \begin{equation*}
    \gamma_{L}(E) = \abs{\Sigma^{\ast} / \Theta(E,L)}
    \stackrel{\ref{lem:entropy-of-automaton-as-entropy-of-cover}(2)}{\leq}
    \abs{X/\Lambda_{E}}
    \stackrel{\ref{lem:entropy-of-automaton-as-entropy-of-cover}(1)}{=}
    \abs{\mathcal{V} \setminus \set{\emptyset}}
    = N(\mathcal{V}) = (E:\mathcal{U})_{\alpha} .
  \end{equation*}
  Finally, if $\mathcal{A}$ is trim, then
  \Cref{lem:entropy-of-automaton-as-entropy-of-cover}~(3) asserts
  $\abs{\Sigma^{\ast} / \Theta(E,L)} = \abs{X/\Lambda_{E}}$ and therefore
  $\gamma_{L}(E) = (E:\mathcal{U})_{\alpha}$.
\end{Proof}

The particular choice of the cover $\mathcal{U} = \set{ F, X \setminus F}$ seems arbitrary, but
this is not the case. Indeed, if the automaton $\mathcal{A} = (Q, \Sigma, \alpha, x_{0}, F)$ is
minimal, then the entropy $\eta(\alpha,\Sigma \cup \set{\epsilon})$ of the automaton equals
$\eta(\alpha, \Sigma \cup \set{\epsilon}, \mathcal{U})$.  We shall show this fact in
\Cref{thm:entropy-of-minimal-automaton-given-by-cover-U}.  As a preparation, we shall first
investigate three auxiliary statements.

\begin{Lemma}
  \label{lem:composition-of-covers}
  Let $X$ be a set, let $S$ be a semigroup, and let $\alpha \colon S \times X \to X$ be an
  action of $S$ on $X$. Let $\mathcal{U}$ be a finite cover of $X$ and let
  $M, N \subseteq S$ be finite. Then
  \begin{equation*}
    \bigvee_{s \in MN} s^{-1}(\mathcal{U}) \equiv \bigvee_{s \in N}
    s^{-1} \bigl(\bigvee_{t \in M} t^{-1}(\mathcal{U})\bigr).
  \end{equation*}
  In particular, the complexities of those two covers coincide.
\end{Lemma}

\begin{Proof}
  Without loss of generality we may assume that $\mathcal{U}$ is closed under
  intersection: in fact, $\mathcal{U}$ is refinement-equivalent to the finite cover
  $\tilde{\mathcal{U}} \coloneqq \set{ \bigcap \mathcal{V} \mid \mathcal{V} \subseteq
    \mathcal{U}}$.
  Hence, if the desired statement was true for $\tilde{\mathcal{U}}$, then this would
  imply
  \begin{equation*}
    \bigvee_{s \in MN} s^{-1}(\mathcal{U}) \equiv \bigvee_{s \in MN}
    s^{-1}(\tilde{\mathcal{U}}) \equiv \bigvee_{s \in N} s^{-1} \bigl(\bigvee_{t \in M}
    t^{-1}(\tilde{\mathcal{U}})\bigr) \equiv \bigvee_{s \in N} s^{-1} \bigl(\bigvee_{t \in M}
    t^{-1}(\mathcal{U})\bigr)
  \end{equation*}
  due to the statements (2) and (3) of Remark~\ref{rem:entropy}.

  Henceforth, assume that $\mathcal{U}$ is closed under intersection. We shall show an
  even stronger claim, namely
  \begin{equation}
    \label{eq:5}
    \bigvee_{s \in MN} s^{-1}(\mathcal{U})
    = \bigvee_{s \in N} s^{-1} \bigl(\bigvee_{t \in M} t^{-1}(\mathcal{U})\bigr).
  \end{equation}
  To ease readability, let us denote the left-hand side by $\mathcal{L}$, and the
  right-hand side by $\mathcal{R}$.

  Let $Y \in \mathcal{L}$.  Then
  \begin{equation*}
    Y = \bigcap_{s \in MN} s^{-1}(U_{s})
  \end{equation*}
  for some $(U_{s} \mid s \in MN) \in \prod_{s \in MN} s^{-1}(\mathcal{U})$.  For each
  $s \in MN$ we can choose $\tau_{s} \in M, \sigma_{s} \in N$ such that
  $s = \tau_{s}\sigma_{s}$.  Then
  \begin{align*}
    Y
    &= \bigcap_{s \in MN} s^{-1}(U_{s}) \\
    &= \bigcap_{s \in MN} (\tau_{s}\sigma_{s})^{-1}(U_{s}) \\
    &= \bigcap_{s \in MN} \sigma_{s}^{-1}\bigl(\tau_{s}^{-1}(U_{\tau_{s}\sigma_{s}})\bigr) \\
    &= \bigcap_{\sigma \in N} \bigcap_{\tau \in M} \sigma^{-1}\bigl( \tau^{-1}(U_{\tau\sigma})
      \bigr) \\
    &= \bigcap_{\sigma \in N} \sigma^{-1} \bigl( \bigcap_{\tau \in M}
      \tau^{-1}(U_{\tau\sigma})\bigr) \in \mathcal{R}
  \end{align*}

  Conversely, let $Y \in \mathcal{R}$.  Then
  \begin{equation*}
    Y = \bigcap_{\sigma \in N} \sigma^{-1}\bigl( \bigcap_{\tau \in M} \tau^{-1}( U_{\sigma, \tau})
    \bigr)
  \end{equation*}
  for some $(U_{\sigma, \tau} \mid \sigma \in N, \tau \in M) \in \prod_{(\sigma,\tau) \in M\times N}
  U^{M \times N}$.  Then
  \begin{align*}
    Y
    &= \bigcap_{\sigma \in N} \bigcap_{\tau \in M}
      \sigma^{-1}\bigl(\tau^{-1}(U_{\sigma,\tau})\bigr)\\
    &= \bigcap_{\sigma \in N} \bigcap_{\tau \in M} (\tau\sigma)^{-1} (U_{\sigma, \tau}) \\
    &= \bigcap_{s \in MN} s^{-1} \Bigl( \bigcap \bigl\{ U_{\sigma, \tau} \mid \sigma \in N, \tau \in
      M, s = \tau\sigma \bigr\}\Bigr)
  \end{align*}
  Define
  \begin{equation*}
    U_{s} := \bigcap \bigl\{ U_{\sigma, \tau} \mid \sigma \in N, \tau \in M, s = \tau\sigma \bigr\}.
  \end{equation*}
  Then $U_{s} \in \mathcal{U}$ for each $s \in MN$, as $\mathcal{U}$ is closed under
  intersections. But then
  \begin{equation*}
    Y = \bigcap_{s \in MN} s^{-1}(U_{s}) \in \mathcal{L}
  \end{equation*}
  as required.

  Finally, Equation~\ref{eq:5} and Remark~\ref{rem:entropy}~(1) yield
  \begin{equation*}
    N\bigl(\bigvee_{s \in MN} s^{-1}(\mathcal{U})\bigr)
    = N\Bigl(\bigvee_{s \in N} s^{-1} \bigl(\bigvee_{t \in M} t^{-1}(\mathcal{U})\bigr)\Bigr),
  \end{equation*}
  as it has been claimed.
\end{Proof}

\begin{Lemma}
  \label{lem:refinement}
  Let $L \subseteq \Sigma^{*}$ and let $\mathcal{A} = (X, \Sigma,\alpha, x_{0},F)$ be the
  minimal automaton of $L$. Consider $S \coloneqq \Sigma \cup \{ \epsilon \}$ and
  $\mathcal{U} \coloneqq \{ F, \, X\setminus F \}$. If $\mathcal{V}$ is a finite open
  cover of $X$, then there exists some $n \in \NN$ such that
  $\mathcal{V} \preceq \bigvee_{s \in S^{n}} s^{-1}(\mathcal{U})$.
\end{Lemma}

\begin{Proof}
  For $n \in \NN$, let us consider the equivalence relation
  \begin{equation*}
    \Lambda_{n} \coloneqq \Lambda_{\Sigma^{(n)}} = \set{ (x, y) \in X \times X \mid \forall w \in
      \Sigma^{(n)} \colon \alpha(w,x) \in F \iff \alpha(w, y) \in F }
  \end{equation*}
  (cf.~Lemma~\ref{lem:entropy-of-automaton-as-entropy-of-cover}).  We are going to show
  that
  \begin{equation*}
    \mathcal{W} \coloneqq \{ [x]_{\Lambda_{n}} \mid n \in \NN, \, x \in X, \, \exists V \in
    \mathcal{V} \colon \, [x]_{\Lambda_{n}} \subseteq V \}
  \end{equation*}
  is an open cover of $X$. By
  Lemma~\ref{lem:entropy-of-automaton-as-entropy-of-cover}~(1), it follows that
  $\mathcal{W}$ is a collection of open subsets of $X$. Thus, we only need to argue that
  $X = \bigcup \mathcal{W}$. To this end, let $x \in X$. Since $\mathcal{V}$ is a cover of
  $X$, there exists some $V \in \mathcal{V}$ with $x \in V$. As $V$ is open in $X$ with
  respect to the subspace topology inherited from $\{ 0,1 \}^{\Sigma^{\ast}}$, we find a
  finite set $E \subseteq \Sigma^{\ast}$ such that
  $W \coloneqq \{ y \in X \mid \forall w \in E \colon \, x(w) = y(w) \} \subseteq V$. Let
  $n \in \NN$ where $E \subseteq S^{n}$. We observe that
  \begin{align*}
    [x]_{\Lambda_{n}}
    &= \{ y \in X \mid \forall w \in S^{n} \colon \, \alpha (x,w) \in F \iff \alpha (y,w) \in F \}
  \\
    &= \{ y \in X \mid \forall w \in S^{n} \colon \, \alpha (x,w)(\epsilon) = 1 \iff \alpha
      (y,w)(\epsilon) = 1 \} \\
    &= \{ y \in X \mid \forall w \in S^{n} \colon \, x(w) = 1 \iff y(w) = 1 \} \\
    &= \{ y \in X \mid \forall w \in S^{n} \colon \, x(w) = y(w) \} \\
    &\subseteq W \subseteq V .
  \end{align*}
  Accordingly, $[x]_{\Lambda_{n}} \in \mathcal{W}$ and hence $x \in \bigcup
  \mathcal{W}$.
  This proves the claim. Now, since $X$ is compact, there exists a finite subset
  $\mathcal{W}_{0}$ where $X = \bigcup \mathcal{W}_{0}$. Due to finiteness of
  $\mathcal{W}_{0}$, there is some $n \in \NN$ such that
  $\mathcal{W}_{0} \preceq X/\Lambda_{n}$. We conclude that
  \begin{equation*}
    \mathcal{V} \preceq \mathcal{W} \preceq \mathcal{W}_{0} \preceq X/\Lambda_{n}
    \stackrel{\ref{lem:entropy-of-automaton-as-entropy-of-cover}(1)}{=} \left( \bigvee\nolimits_{s
        \in S^{n}} s^{-1}(\mathcal{U}) \right) \setminus \{ \emptyset \} \equiv \left(
      \bigvee\nolimits_{s \in S^{n}} s^{-1}(\mathcal{U}) \right) ,
  \end{equation*}
  which completes the proof.
\end{Proof}

We finally reached the point where we can show that the Myhill-Nerode complexity
and the topological entropy of $L$ coincide.

\begin{Theorem}
  \label{thm:entropy-of-minimal-automaton-given-by-cover-U}
  Let $L \subseteq \Sigma^{*}$ and let $\mathcal{A} = (X, \Sigma,\alpha, x_{0},F)$ be the minimal
  automaton of $L$. Consider $S \coloneqq \Sigma \cup \{ \epsilon \}$ and
  $\mathcal{U} \coloneqq \{ F, \, X\setminus F \}$. Then
  $h(L) = \eta (\alpha,S,\mathcal{U}) = \eta (\alpha, S)$.
\end{Theorem}

\begin{Proof}
  Define $\mathcal{U} \coloneqq \set{ F, \, X \setminus F }$. Since $\mathcal{A}$ is trim,
  we know that $h(L) = \eta (\alpha,S,\mathcal{U})$ by
  Theorem~\ref{thm:entropy-of-automaton-as-entropy-of-cover} and hence
  $h(L) \leq \eta (\alpha,S)$.  To show the converse inequality, let $\mathcal{V}$ be a
  finite open cover of $X$. We show that
  $\eta (\alpha,S,\mathcal{V}) \leq \eta (\alpha,S,\mathcal{U})$. According to
  \Cref{lem:refinement}, there exists some $m \in \NN$ such that $\mathcal{V}$ is refined
  by $\bigvee_{s \in S^{m}} s^{-1}(\mathcal{U})$. Then
  \begin{equation*}
    N\Bigl( \bigvee\nolimits_{s \in S^{n}} s^{-1}(\mathcal{V}) \Bigr) \leq
    N\Bigl( \bigvee\nolimits_{s \in S^{n}} s^{-1}\bigl( \bigvee\nolimits_{t \in S^{m}}
        t^{-1}(\mathcal{U}) \bigr) \Bigr) = N \Bigl( \bigvee\nolimits_{s \in S^{m+n}}
      s^{-1}(\mathcal{U}) \Bigr)
  \end{equation*}
  by \Cref{lem:composition-of-covers}. Now we obtain
  \begin{align*}
    \eta (\alpha, S, \mathcal{V})
    &= \limsup_{n \to \infty} \frac{\log_{2} (S^{n} : \mathcal{V})_{\alpha}}{n}
      \stackrel{\ref{rem:entropy}(1)}{\leq} \limsup_{n \to \infty} \frac{ \log_{2} (S^{n+m} :
      \mathcal{U})_{\alpha}}{n} \\
    &= \limsup_{n \to
      \infty} \frac{ \log_{2} (S^{n} : \mathcal{U})_{\alpha}}{n} = \eta (\alpha, S, \mathcal{U}).
  \end{align*}
  Therefore, $\eta (\alpha ,S) \leq \eta(\alpha, S, \mathcal{U})$ and hence
  $\eta (\alpha,S) = \eta (\alpha,S,\mathcal{U}) = h(L)$ by
  \Cref{thm:entropy-of-automaton-as-entropy-of-cover}.
\end{Proof}

\section{Examples}
\label{sec:examples}

\Cref{thm:entropy-of-minimal-automaton-given-by-cover-U} allows us to easily
compute the topological entropy of certain classes of languages.  To begin with,
we show that all regular languages have zero entropy.

\begin{Theorem}
  Let $\Sigma$ be an alphabet and $L \subseteq \Sigma^{\ast}$. The following are
  equivalent:
  \begin{enumerate}
  \item\label{item:1} $L$ is regular,
  \item\label{item:2} $\gamma_{L}$ is bounded, and
  \item\label{item:3} there exists some finite subset $F \subseteq
    \Sigma^{\ast}$ such that $\Theta(F,L) = \Theta(L)$.
  \end{enumerate}
\end{Theorem}

\begin{Proof}
  $\text{\ref{item:1}}\Longrightarrow\text{\ref{item:2}}$. Due to \Cref{Theorem:Myhill-Nerode},
  $\Theta (L)$ has finite index. Note that $\Theta (L) \subseteq \Theta (F,L)$ and hence
  $\gamma_{L}(F) \leq \ind \Theta (L)$ for all $F \subseteq \Sigma^{\ast}$ finite. Thus,
  $\gamma_{L}$ is bounded.

  $\text{\ref{item:2}}\Longrightarrow\text{\ref{item:3}}$. Suppose that $\gamma_{L}$ is bounded. Then
  there exists some finite $F_{0} \subseteq \Sigma^{\ast}$ such that
  $\gamma_{L}(F_{0}) = \sup \set{ \gamma_{L}(F) \mid F \subseteq \Sigma^{\ast} \text{
      finite}}$.
  We shall show that $\Theta(F_{0},L) = \Theta(L)$. Of course,
  $\Theta (L) \subseteq \Theta (F_{0},L)$.  Let
  $(u,v) \in (\Sigma^{\ast} \times \Sigma^{\ast}) \setminus \Theta (L)$. By~\eqref{eq:1}
  there exists some finite $F_{1} \subseteq \Sigma^{\ast}$ such that
  $(u,v) \notin \Theta (F_{1},L)$.  Obviously, $F_{0} \cup F_{1} \subseteq \Sigma^{\ast}$
  is finite and $\Theta (F_{0} \cup F_{1},L) \subseteq \Theta (F_{0},L)$. By assumption,
  $\gamma_{L}(F_{0} \cup F_{1}) \leq \gamma_{L}(F_{0})$.  Consequently,
  $\Theta (F_{0} \cup F_{1},L) = \Theta (F_{0},L)$ and therefore
  $(u,v) \in (\Sigma^{\ast} \times \Sigma^{\ast}) \setminus \Theta (F_{1},L) \subseteq
  (\Sigma^{\ast} \times \Sigma^{\ast}) \setminus \Theta (F_{0} \cup F_{1},L) =
  (\Sigma^{\ast} \times \Sigma^{\ast}) \setminus \Theta (F_{0},L)$.
  This substantiates that $\Theta(F_{0},L) = \Theta(L)$.

  $\text{\ref{item:3}}\Longrightarrow\text{\ref{item:1}}$. By assumption $\Theta(L) = \Theta(F, L)$, and
  since $\Theta(F, L)$ has finite index, $\Theta (L)$ has finite index as well.  Hence,
  $L$ is regular due to \Cref{Theorem:Myhill-Nerode}.
\end{Proof}

\begin{Corollary}
  Let $\Sigma$ be an alphabet.  If $L \subseteq \Sigma^{\ast}$ is regular, then $h(L) = 0$.
\end{Corollary}

The converse of this corollary does not hold, i.e., there are non-regular languages with
vanishing topological entropy.  To see this we shall show that \emph{Dyck languages}
always have zero entropy (cf.\,~\Cref{corollary:dyck.language}).  We shall put the
corresponding argumentation in a more general framework, by estimating the entropy of
languages defined by groups. For this purpose, we recall the concept of \emph{growth} in
groups. Consider a finitely generated group $G$. Let $S$ be a finite symmetric generating
subset of $G$ containing the neutral element. The \emph{exponential growth rate} of $G$
with respect to $S$ is defined to be
\begin{displaymath}
  \egr (G,S) \coloneqq \limsup_{n \to \infty} \frac{\log_{2} \vert S^{n} \vert}{n}.
\end{displaymath}
Note that this quantity is finite as $\vert S^{n} \vert \leq \vert S \vert^{n}$ for every
$n \in \mathbb{N}$. Furthermore,
\begin{displaymath}
  \egr (G,S) = \lim_{n \to \infty} \frac{\log_{2} \vert S^{n} \vert}{n}
\end{displaymath}
due to a well-known result by Fekete~\cite{Fekete}. Of course, the precise value of the exponential
growth rate depends upon the particular choice of a generating set.


However, if $T$ is another finite symmetric generating subset of $G$ containing the
neutral element, then
\begin{displaymath}
  \frac{1}{k} \cdot \egr (G,T) \leq \egr (G,S) \leq l \cdot \egr (G,T)
\end{displaymath}
where $k \coloneqq \inf \{ m \in \mathbb{N}\setminus \{ 0 \} \mid T \subseteq S^{m} \}$ and
$l \coloneqq \inf \{ m \in \mathbb{N}\setminus \{ 0 \} \mid S \subseteq T^{m} \}$. This justifies
the following definition: $G$ is said to have \emph{sub-exponential growth} if $\egr (G,S) = 0$ for
some (and thus any) symmetric generating set $S$ of $G$ containing the neutral element. The class of
finitely generated groups with sub-exponential growth encompasses all finitely generated abelian
groups. In fact, if $G$ is abelian, then
\begin{equation*}
  S^{n} \subseteq \Bigl\{ \prod_{s \in S} s^{\alpha(s)} \Bigm| \alpha \colon S \to \set{0, \dots, n} \Bigr\}
\end{equation*}
and thus $\abs{S^{n}} \leq (n+1)^{\abs{S}}$ for all $n \in \mathbb{N}$. Now let us return to formal languages.

\begin{Theorem}%
  \label{theorem:generalized.dyck.language}%
  Let $\Sigma$ be an alphabet. Let $G$ be a group, $\phi \colon \Sigma^{\ast} \to G$ a
  homomorphism, $H \subseteq G$, and $E \subseteq G$ finite. Define
  \begin{align*}
    P_{\phi}(H) &\coloneqq \{ w \in \Sigma^{\ast} \mid \forall u \textit{ prefix of } w
                  \colon\phi(u) \in H \} , \\
    L_{\phi}(H,E) &\coloneqq P_{\phi}(H) \cap \phi^{-1}(E).
  \end{align*}
  Then $\gamma (F,L_{\phi}(H,E)) \leq \abs{E}\cdot \abs{\phi(F)} + 1$ for all finite
  $F \subseteq \Sigma^{\ast}$.  In particular,
  \begin{equation*}
    h(L_{\phi}(H,E)) \leq \limsup_{n \to \infty} \frac{\log_{2}|\phi (\Sigma^{\left( n
        \right)})|}{n} \leq \log_{2} |\Sigma | .
  \end{equation*} Furthermore, if $S$ is a finite symmetric generating subset of $G$ containing the neutral element and $k \coloneqq \inf \{ m \in \mathbb{N}\setminus \{ 0 \} \mid \phi (\Sigma) \subseteq S^{m} \}$, then \begin{equation*}
  h(L_{\phi}(H,E)) \leq k \cdot \egr (G,S) .
  \end{equation*}
\end{Theorem}
\begin{Proof}
  We abbreviate $P \coloneqq P_{\phi}(H)$ and $L \coloneqq L_{\phi}(H,E)$. Consider a finite subset
  $F \subseteq \Sigma^{\ast}$. Then $Q \coloneqq E\phi(F)^{-1}$ is a finite subset of $G$. Fix any
  object $\infty \notin Q$ and define $Q_{\infty} \coloneqq Q \cup \{ \infty \}$.  Let us consider
  the map $\psi \colon \Sigma^{\ast} \to Q_{\infty}$ given by
  \begin{equation*}
    \psi (u) \coloneqq \begin{cases}
      \phi (u) & \textnormal{if } u \in P \cap \phi^{-1}(Q) , \\
      \infty & \textnormal{otherwise}
    \end{cases} \qquad (u \in \Sigma^{\ast}) .
  \end{equation*}
  We show $\ker \psi \subseteq \Theta(F,L)$. To this end, let
  $(u,v) \in \ker \psi$. We proceed by case analysis.

  First case: $\psi (u) = \psi (v) \ne \infty$. Now, $u,v \in P \cap \phi^{-1}(Q)$ and
  $\phi(u) = \psi(u) = \psi(v) = \phi(v)$. Let $w \in F$ and suppose that $uw \in L$. We show
  $vw \in L$. We observe that
  \begin{equation*}
    \phi(vw) = \phi(v)\phi(w) = \phi(u)\phi(w) = \phi(uw) \in E,
  \end{equation*}
  i.e., $vw \in \phi^{-1}(E)$. In order to prove that $vw \in P$, let $x$ be a prefix of
  $vw$. If $x$ is a prefix of $v$, then $\phi(x) \in H$ as $v \in P$. Otherwise, there
  exists a prefix $y$ of $w$ such that $x = vy$, and so we conclude that
  $\phi(x) = \phi(vy) = \phi(v)\phi(y) = \phi(u)\phi(y) = \phi(uy) \in H$, because
  $uw \in P$ and $uy$ is a prefix of $uw$. Hence, $vw \in L$. On account of symmetry, it
  follows that $(u,v) \in \Theta(F,L)$.

  Second case: $\psi (u) = \psi (v) = \infty$. Let $x \in \{ u, \, v \}$. If
  $x \notin \phi^{-1}(Q)$, then we conclude that $\phi (xw) = \phi (x)\phi (w) \notin E$ and thus
  $xw \notin L$ for any $w \in F$. If $x \notin P$, then $xw \notin P$ and hence
  $xw \notin L$ for any $w \in F$. This proves that
  $\{ uw, \, vw \} \cap L = \emptyset$ for all $w \in F$. Consequently,
  $(u,v) \in \Theta (F,L)$.

  This substantiates that $\ker \psi \subseteq \Theta(F,L)$.  Therefore
  \begin{equation*}
    \gamma (F,L) = \ind \Theta (F,L) \leq \ind (\ker \psi) \leq |Q_{\infty}|
    \leq |Q| + 1 \leq |E|\cdot |\phi (F)| + 1 .
  \end{equation*}
  In particular, it follows that
  \begin{align*}
    h(L)
    &= \limsup_{n \to \infty}\frac{\log_{2}\gamma_{L} (\Sigma^{\left( n \right)})}{n}
      \leq \limsup_{n \to \infty}\frac{\log_{2} (|E|\cdot |\phi (\Sigma^{\left( n \right)})| + 1)}{n} \\
    &= \limsup_{n \to \infty}\frac{\log_{2} |\phi (\Sigma^{\left( n \right)})|}{n}
      \leq \limsup_{n \to \infty}\frac{\log_{2} (|\Sigma |^{n})}{n} = \log_{2} |\Sigma | .
  \end{align*}
  Finally, suppose $S$ to be a finite symmetric generating subset of $G$ containing the neutral
  element. Since $\Sigma$ is finite,
  $M \coloneqq \{ m \in \mathbb{N}\setminus \{ 0 \} \mid \phi(\Sigma) \subseteq S^{m}\}$ is not
  empty. Let $k \coloneqq \inf M$. Our considerations above now readily imply that
  \begin{equation*}
    h(L_{\phi}(S, E)) \leq \limsup_{n \to \infty} \frac{\log_{2}\lvert \phi(\Sigma^{(n)}) \rvert}{n}
    \leq k \cdot \limsup_{n \to \infty}\frac{\log_{2}\lvert S^{n} \rvert}{n} = k \cdot \egr (G,S).
  \end{equation*}
\end{Proof}

For groups whose growth is sub-exponential the previous theorem yields that the corresponding
languages $L_{\phi}(S, E)$ have zero entropy.

\begin{Corollary}
  \label{corollary:zero-entropy-for-subexponential-group-growth}%
  Let $\Sigma$ be an alphabet, let $G$ be a group with sub-exponential growth, and
  $\phi\colon \Sigma^{*} \to G$ a homomorphism.  Then for each $S \subseteq G$ and finite
  $E \subseteq G$, it is true that $h(L_{\phi}(S, E)) = 0$.
\end{Corollary}

We immediately obtain the following statement.

\begin{Corollary}\label{corollary:abelian.dyck.language}%
  Let $\Sigma$ be an alphabet, let $G$ be a finitely generated abelian group, and
  $\phi \colon \Sigma^{\ast} \to G$ a homomorphism.  Then for each $S \subseteq G$ and finite
  $E \subseteq G$, it is true that $h(L_{\phi}(S,E)) = 0$.
\end{Corollary}

The following corollaries are immediate consequences of
Theorem~\ref{theorem:generalized.dyck.language} for $S = G$.

\begin{Corollary}%
  \label{corollary:group}%
  Let $\Sigma$ be a finite alphabet and $L \subseteq \Sigma^{\ast}$. Let $G$ be a group,
  $\phi \colon \Sigma^{\ast} \to G$ a homomorphism and $E \subseteq G$ finite such that
  $L = \phi^{-1}(E)$. Then $\gamma (F,L) \leq \abs{E} \cdot \abs{\phi(F)} + 1$ for all finite
  $F \subseteq \Sigma^{\ast}$. In particular,
  \begin{equation*}
    h(L) \leq \limsup_{n \to \infty} \frac{\log_{2}|\phi (\Sigma^{\left( n
        \right)})|}{n} \leq \log_{2} |\Sigma | .
  \end{equation*}
\end{Corollary}

\begin{Corollary}\label{corollary:abelian.group}%
  Let $\Sigma$ be a finite alphabet, $L \subseteq \Sigma^{\ast}$. Let $G$ be an abelian group,
  $\phi \colon \Sigma^{\ast} \to G$ a homomorphism and $E \subseteq G$ finite such that
  $L = \phi^{-1}(E)$. Then $h(L) = 0$.
\end{Corollary}

With the previous results in place, we are now able to argue that \emph{Dyck languages}
have finite entropy.  Recall that the \emph{Dyck language with $k$ sorts of parentheses}
consists of all balanced strings over $\set{ (_{1}, )_{1}, \dots, (_{k}, )_{k} }$.
Alternatively, we can view the Dyck language with $k$ sorts of parentheses as the set of
all strings that can be reduced to the empty word by successively eliminating matching
pairs of parentheses.

We can formalize this as follows.  Let $\Sigma, \overline{\Sigma}$ be two alphabets,
$\Delta \coloneqq \Sigma \cup \overline{\Sigma}$, and let
$\kappa \colon \Sigma \to \overline{\Sigma}$ be a bijection.  Consider the the free group
$F(\Sigma)$ with generator set $\Sigma$, and denote with
$\phi \colon \Delta^{*} \to F(\Sigma)$ the unique homomorphism satisfying $\phi(a) = a$ and
$\phi(\kappa(a)) = a^{-1}$ for all $a \in \Sigma$.  Define
\begin{equation*}
  D(\kappa) \coloneqq \{ w \in \Delta^{\ast} \mid
  \phi (w) = e \wedge (\forall u \text{ prefix of } w: |u|_{a} \geq |u|_{\kappa(a)}) \}.
\end{equation*}
If $\Sigma = \set{ (_{1}, \dots, (_{k}}$, $\overline{\Sigma} = \set{)_{1}, \dots, )_{k}}$,
and $\kappa\bigl((_{i}\bigr) = {)_{i}}$, then the set $D(\kappa)$ coincides with the Dyck
language with $k$ sorts of parentheses.

\begin{Theorem}\label{corollary:dyck.language}%
  Let $\kappa \colon \Sigma \to \overline{\Sigma}$ be a bijection between finite sets.
  Then
  \begin{equation*}
    \log_{2}\abs{\Sigma} \leq h(D(\kappa)) \leq \log_{2}(2\lvert \Sigma \rvert - 1)
  \end{equation*}
  for $S \coloneqq \Sigma \cup \Sigma^{-1} \cup \{ e \}$, where $e$ denotes the neutral
  element of $F(\Sigma)$.
\end{Theorem}
\begin{Proof}%
  Let $L \coloneqq D(\kappa)$.  We first show the inequality $\ind \Theta(\Sigma^{(n)}, L)
  \geq \abs{\Sigma^{n}}$, since this implies $\log_{2}\abs{\Sigma} \leq h(D(\kappa))$.
  For this let $u, v \in \Sigma^{n}$, $u \neq v$.  Define $\kappa(u) :=
  \kappa(u_{\abs{u}}) \dots \kappa(u_{1})$, where $u = u_{1} \dots u_{\abs{u}}$.  Then $u
  \cdot \kappa(u) \in L$, but $v \cdot \kappa(u) \notin L$.  Thus $(u, v) \notin
  \Theta(\Sigma^{(n)}, L)$ and therefore $\ind \Theta(\Sigma^{(n)}, L) \geq
  \abs{\Sigma^{n}}$ as required.

  For the second inequality let us consider the unique homomorphism
  $\psi \colon F(\Sigma) \to \mathbb{Z}^{\Sigma}$ satisfying
  \begin{equation*}
    \psi(b)(a) \coloneqq \begin{cases}
      1 & \text{if } a= b, \\
      0 & \text{otherwise}
    \end{cases}
    \qquad (a, b \in \Sigma) .
  \end{equation*}
  We observe $D(\kappa) = L_{\phi}(\psi^{-1}(\mathbb{N}^{\Sigma}),\{ e \})$, where the
  mapping $\phi$ is as above.  Hence, we have $h(D(\kappa)) \leq \egr (F(\Sigma),S)$ by
  \Cref{corollary:abelian.dyck.language}.  As it is known that $\egr(F(\Sigma), S) =
  \log_{2}(2\lvert \Sigma \vert - 1)$ we obtain the claim.
\end{Proof}

Note that for $\abs{\Sigma} = 1$ we have $h(D(\kappa)) = 0$.  Thus $D(\kappa)$ is an
example of a non-regular language with zero entropy.  For $\abs{\Sigma} > 1$ the exact
value of $h(D(\kappa))$ is unknown to the authors.

The reason that Dyck languages with more than one type of parentheses have non-zero
positive entropy is the following: the different types of parentheses occurring in a word
$w \in D(\kappa)$ need to be mutually balanced, i.e., $\phi(w) = e$.  In other words, if
we replace this requirement by the weaker condition that each opening parenthesis has to
be closed eventually, then we obtain a class of languages with zero entropy.

\begin{Theorem}\label{corollary:fake.dyck.language}%
  Let $\kappa \colon \Sigma \to \overline{\Sigma}$ be a bijection between finite sets, let
  $\Delta \coloneqq \Sigma \cup \overline{\Sigma}$, and consider the language
  \begin{equation*}
    D'(\kappa) \coloneqq \{ w \in \Delta^{\ast} \mid \forall a \in \Sigma \colon (|w|_{a} =
    |w|_{\kappa(a)}) \wedge (\forall u \textit{ prefix of } w: |u|_{a} \geq |u|_{\kappa(a)}) \} .
  \end{equation*}
  Then $h(D'(\kappa)) = 0$.
\end{Theorem}
\begin{Proof}%
  Let us consider the homomorphism $\phi \colon \Delta^{\ast} \to \mathbb{Z}^{\Sigma}$
  given by
  \begin{equation*}
    \phi(w)(a) \coloneqq |w|_{a} - |w|_{\kappa(a)} \qquad (w \in \Delta^{\ast}, \, a \in \Sigma ) .
  \end{equation*}
  We observe that $D(\kappa) = L_{\phi}(\mathbb{N}^{\Sigma},\{ 0 \})$, wherefore
  $h(D(\kappa)) = 0$ by \Cref{theorem:generalized.dyck.language}.
\end{Proof}

Other non-regular languages with vanishing entropy are discussed in the following
examples.

\begin{Example} Let $\Sigma$ be an alphabet.
  \begin{enumerate}
  \item Let $m \in \mathbb{N}$ and $a,b \in \Sigma$, $a \ne b$. Then $L
    \coloneqq \{ w \in \Sigma^{\ast} \mid \vert w \vert_{a} = \vert w \vert_{b}
    + m \}$ is not regular. However, $h(L) = 0$ by
    Corollary~\ref{corollary:abelian.group}. To see this, note that the mapping
    $\phi \colon \Sigma^{\ast} \to \mathbb{Z}, \, w \mapsto |w|_{a} - |w|_{b}$
    constitutes a homomorphism where $L = \phi^{-1}(\{ m \})$.

  \item Suppose $\Sigma = \set{ a,b,c }$. Then $L \coloneqq \set{
      a^{m}b^{m}c^{m} \mid m \in \mathbb{N} }$ is not context-free, but $h(L) =
    0$.  To see this we show that for every $n$ the relation $\Theta =
    \Theta(\Sigma^{(n)}, L)$ has the equivalence classes
    \begin{equation}
      \label{eq:2}
      \begin{split}
        [a^{k}]_{\Theta}, & \quad k \leq n/2\\
        [a^{k}b^{\ell}]_{\Theta}, & \quad 1 \leq \ell \leq k,\, 2k - \ell \leq n\\
        [a^{k}b^{k}c^{\ell}]_{\Theta}, & \quad 1 \leq \ell \leq k,\, k - \ell \leq n\\
        [b]_{\Theta}.
      \end{split}
    \end{equation}
    From this it follows $\ind \Theta(\Sigma^{(n)}, L) \in \mathcal{O}(n^{2})$,
    and thus $h(L) = 0$.

    To see that the sets in~\eqref{eq:2} are indeed all equivalence classes of
    $\Theta(\Sigma^{(n)}, L)$, let $u \in \Sigma^{*}$ such that $u$ is not an
    element of the first three types of classes in~\eqref{eq:2}.  We need to
    show that then $u \in [b]_{\Theta(\Sigma^{(n)}, L)}$.  We do this by showing
    that there is no $w \in \Sigma^{(n)}$ such that $uw \in L$.

    Assume by contradiction that such a word $w$ exists.  Then $w$ must be of
    one of the following forms
    \begin{align*}
      w &= a^{\ell}b^{k}c^{k}, \quad 2k+\ell \leq n,\, 0 \leq \ell \leq k,\\
      w &= b^{\ell}c^{k}, \quad k + \ell \leq n,\, 0 \leq \ell < k,\\
      w &= c^{\ell}, \quad 0 \leq \ell < n
    \end{align*}

    If $w = a^{\ell}b^{k}c^{k}$, $2k + \ell \leq n$, $0 \leq \ell \leq k$, then
    $u = a^{k-\ell}$, $k-\ell \leq n/2$, and therefore $u \in
    [a^{k-\ell}]_{\Theta(\Sigma^{(n)}, L)}$, a contradiction.  If $w =
    b^{\ell}c^{k}$, $k + \ell \leq n, 0 \leq \ell < k$, then $u =
    a^{k}b^{k-\ell}$, and $k - \ell > 0$, $2k - (k-\ell) \leq n$, thus $u \in
    [a^{k}b^{k-\ell}]_{\Theta(\Sigma^{(n)}, L)}$, again a contradiction.  If $w
    = c^{\ell}$, then $u = a^{k}b^{k}c^{k-\ell}$, and $k - (k - \ell) \leq n$,
    so $u \in [a^{k}b^{k}c^{\ell}]_{\Theta(\Sigma^{(n)}, L)}$, another
    contradiction.

    Thus, our assumption that $w$ exists is false.  The same is true for the
    word $b$, and thus $u \in [b]_{\Theta(\Sigma^{(n)}, L)}$, as required.
  \end{enumerate}
\end{Example}

Next we are looking for an example of a language with non-zero entropy.  Of
course, by what we have already shown, a suitable candidate for this has to be
non-regular.  But do not have to require much more: the following example shows
that there are context-free languages with non-zero entropy.

\begin{Example}
  \label{expl:palindromes}%
  Suppose $\vert \Sigma \vert \geq 2$. Then the \emph{palindrome language}
  \begin{equation*}
    L \coloneqq \{ ww^{R} \mid w \in \Sigma^{\ast} \}
  \end{equation*}
  is not regular, but context-free, and $h(L) \in (0,\infty)$.

  To see $h(L) > 0$, observe that for each $n \in \NN$ and all $u, v \in \Sigma^{n}$, if
  $(u, v) \in \Theta(\Sigma^{(n)}, L)$, then $u = v$.  This is because if $vv^{R} \in L$,
  we also have $uv^{R} \in L$, and hence $u = v$.  Thus  
  \begin{equation*}
    [u]_{\Theta(\Sigma^{(n)}, L)} \neq [v]_{\Theta(\Sigma^{(n)}, L)} \quad (u \neq v)
  \end{equation*}
  Thus $\ind \Theta(\Sigma^{(n)}, L) \geq \abs{\Sigma^{n}} = \abs{\Sigma}^{n}$, and we
  obtain
  \begin{equation*}
    h(L)
    = \limsup_{n\to\infty} \frac{\log_{2}\abs{\Sigma}^{n}}{n}
    = \log_{2}\abs{\Sigma} > 0.
  \end{equation*}

  To see $h(L) < \infty$ we shall consider the relation $\Theta^{*}$ defined by
  \begin{equation*}
    (u, v) \in \Theta^{*} \iff (u, v) \in \Theta(\Sigma^{(n)}, L) \text{ and } (\abs{u} \leq n
    \iff \abs{v} \leq n).
  \end{equation*}
  Then $\ind \Theta(\Sigma^{(n)}, L) \leq \ind \Theta^{*}$.  We shall show
  \begin{equation*}
    \limsup_{n \to \infty} \frac{\log_{2}(\ind \Theta^{*})}{n} < \infty.
  \end{equation*}

  There are at most $\abs{\Sigma^{(n)}}$ many equivalence classes $[u]_{\Theta^{*}}$ for
  $u \in \Sigma^{*}$, $\abs{u} < n$.  To count the number of equivalence classes for
  $\abs{u} \geq n$ we define
  \begin{equation*}
    \ell_{n}(u) := \set{a_{1}\ldots a_{i} \mid 1 \leq i \leq n, a_{1}, \ldots, a_{i} \in \Sigma, u
      = a_{1}\ldots a_{i} u', u' \in L}.
  \end{equation*}
  Then for $u, v \in \Sigma^{*} \setminus \Sigma^{(n)}$ we have
  \begin{equation}
    \label{eq:4}
    (u, v) \in \Theta^{*} \iff (u, v) \in \Theta(\Sigma^{(n)}, L) \iff \ell_{n}(u) = \ell_{n}(v).
  \end{equation}
  The first equivalence is clear.  To see the second equivalence let $(u, v) \in
  \Theta(\Sigma^{(n)}, L)$, and let $a_{1}\ldots a_{i} \in \ell_{n}(u)$.  By definition of
  $\ell_{n}(u)$ it is then true that $u(a_{1}\ldots a_{i})^{R} \in L$.  Because $(u, v)
  \in \Theta(\Sigma^{(n)}, L)$ we therefore obtain $v(a_{1}\ldots a_{i})^{R} \in L$, i.e.,
  $v$ is of the form $v = a_{1}\ldots a_{i} v'$ for some $v' \in L$.  This yields
  $a_{1}\ldots a_{i} \in \ell_{n}(v)$.  By symmetry we obtain $\ell_{n}(u) = \ell_{n}(v)$
  as required.

  Conversely, assume $\ell_{n}(u) = \ell_{n}(v)$, and let $w \in \Sigma^{(n)}$ be such
  that $uw \in L$.  Because $\abs{u} \geq n$, there exists $u' \in L$ with
  $uw = w^{R} u' w$.  Then $w^{R} \in \ell_{n}(u) = \ell_{n}(v)$, and therefore
  $v = w^{R} v'$ for some $v' \in L$.  But then $vw \in L$.  By symmetry
  $vw \in L \implies uw \in L$ for each $w \in \Sigma^{(n)}$, and therefore
  $(u, v) \in \Theta(\Sigma^{(n)}, L)$, as required.

  Using the characterization from Equation~\eqref{eq:4} we have
  \begin{equation*}
    \Bigl\lvert \faktor{\Sigma^{*}\setminus \Sigma^{(n)}}{\Theta^{*}} \Bigr\rvert = \bigl\lvert
    \set{ \ell_{n}(u) \mid u \in \Sigma^{*}\setminus \Sigma^{(n)}} \bigr\rvert
  \end{equation*}

  Now every set $\ell_{n}(u)$ with $u = u_{1}\ldots u_{k}$, $k \geq n$, can be represented
  by the prefix $u_{1}\ldots u_{n}$ of $u$ of length $n$ together with a tuple
  $t \in \set{0,1}^{n}$ defined by
  \begin{equation*}
    t_{i} = 1 \iff u_{1}\ldots u_{i} \in \ell_{n}(u).
  \end{equation*}
  Therefore,
  \begin{equation*}
    \Bigl\lvert \faktor{\Sigma^{*}\setminus \Sigma^{(n)}}{\Theta^{*}} \Bigr\rvert = \bigl\lvert
    \set{ \ell_{n}(u) \mid u \in \Sigma^{*}\setminus \Sigma^{(n)}} \bigr\rvert
    \leq \abs{\Sigma}^{n} \cdot 2^{n}.
  \end{equation*}
  This yields
  \begin{equation*}
    \ind \Theta^{*}
    = \Bigl\lvert \faktor{\Sigma^{(n)}}{\Theta^{*}} \Bigr\rvert +
    \Bigl\lvert \faktor{\Sigma^{*}\setminus \Sigma^{(n)}}{\Theta^{*}} \Bigr\rvert
    \leq \abs{\Sigma^{(n)}} + \abs{\Sigma}^{n}\cdot 2^{n},
  \end{equation*}
  and thus
  \begin{equation*}
    \limsup_{n \to \infty} \frac{\log_{2}(\ind \Theta^{*})}{n} \leq \log_{2}(2\abs{\Sigma}) < \infty.
  \end{equation*}
\end{Example}

It is unclear to the authors whether the upper bound obtained in the proof of
\Cref{expl:palindromes} is related to the one in \Cref{corollary:dyck.language}.

It is not hard to see that the entropy of a formal language can very well be infinite.
This is illustrated by the following example.
\begin{Example}
  Let $\abs{\Sigma} \geq 2$, and choose mappings
  $\phi_{n} \colon \Sigma^{2^{n}} \to \subsets{\Sigma^{n}}$ for each $n \in \NN$ such that
  $\abs{\im{\phi_{n}}} = \abs{\Sigma}^{2^{n}} = 2^{2^{n}}$.  Then define a language
  $L \subseteq \Sigma^{*}$ by
  \begin{equation*}
    L \cap \Sigma^{m} \coloneqq
    \begin{cases}
      \set{ uv \mid u \in \Sigma^{2^{n}}, v \in \phi_{n}(u)} & \text{if $m = 2^{n} + n$ for some $n
        \in \NN$}, \\
      \emptyset & \text{otherwise}.
    \end{cases}
  \end{equation*}
  Then $2^{2^{n}} \leq \gamma_{L}(n)$, i.e.,
  \begin{equation}
    \label{eq:3}
    2^{2^{n}} \leq \ind \Theta(\Sigma^{n}, L).
  \end{equation}
  To see this we shall show that each word $\phi_{n}(u)$ defines its own equivalence
  class, i.e., for words $u_{0}, u_{1} \in \Sigma^{2^{n}}$ with
  $\phi_{n}(u_{0}) \neq \phi_{n}(u_{1})$ we have
  $(u_{0}, u_{1}) \notin \Theta(\Sigma^{n}, L)$.  This is because if
  $\phi_{n}(u_{0}) \neq \phi_{n}(u_{1})$ we can assume without loss of generality that
  there exists some word $v \in \phi_{n}(u_{0}) \setminus \phi_{n}(u_{1})$.  By
  definition of $L$ we then have $u_{0}v \in L$, but since $\abs{u_{1}v} = 2^{n} + n$ and
  $v \notin \phi_{n}(u_{1})$ we also get $u_{1}v \notin L$.  Thus
  $(u_{0}, u_{1}) \notin \Theta(\Sigma^{n}, L)$.

  But then \eqref{eq:3} implies
  \begin{equation*}
    \limsup_{n \to \infty} \frac{\log_{2} \gamma_{L}(n)}{n}
    \geq \limsup_{n \to \infty} \frac{\log_{2} 2^{2^{n}}}{n} = \infty,
  \end{equation*}
  and thus $h(L) = \infty$.
\end{Example}

\section{Topological entropy and entropic dimension}
\label{sec:topol-entr-entr}

\def\Ball{\operatorname{B}}
\def\edim{\operatorname{dim}}

Another interesting characterization of the entropy of formal languages is in terms of the
entropic dimension of a suitable precompact pseudo-ultrametric space.  For this recall
that a pseudo-metric space $(X, d)$ is called \emph{precompact} if for each $r \in (0,
\infty)$ there exists some finite set $F \subseteq X$ such that
\begin{equation*}
  X = \bigcup\set{ \Ball_{d}(x,r) \mid x \in F }.
\end{equation*}
If $(X, d)$ is a precompact pseudo-metric space, then define
\begin{equation*}
  \gamma_{(X,d)}(r) := \inf\bigl\{\, \abs{F} \bigm| F \subseteq X \text{ finite}, X = \bigcup\set{
    \Ball_{d}(x, r) \mid x \in F } \,\bigr\}.
\end{equation*}
Then the \emph{entropic dimension} $\edim(X, d)$ of the precompact pseudo-metric space
$(X, d)$ is defined as~\cite{journal/jfuncana/Cohen82}
\begin{equation*}
  \edim(X, d) \coloneqq \limsup_{r \to 0+} \frac{\log_{2}(\gamma_{(X,d)}(r))}{\log_{2}(1/r)}.
\end{equation*}

To now obtain a precompact pseudo-metric space $(X, d)$ whose entropic dimension is the
same as the topological entropy of a given language $L$, we shall first start with a
general observation.  Let $X$ be a non-empty set and let
$\Theta = ( \Theta_{n} \mid n \in \NN )$ be a descending sequence of equivalence relations
on $X$.  Define $d_{\Theta} \colon X \times X \to [0, \infty)$ as
\begin{equation*}
  d_{\Theta}(x, y) \coloneqq 2^{-\inf\set{n \in \NN \mid (x, y) \notin \Theta_{n}}} \quad (x, y \in X).
\end{equation*}
It is easy to see that $d_{\Theta}(x, x) = 0$ and $d_{\Theta}(x, y) = d_{\Theta}(y, x)$ is
true for all $x, y \in X$.  Moreover, as
\begin{align*}
  \set{ n \in \NN \mid (x, z) \notin \Theta_{n} }
  &\subseteq \set{ n \in \NN \mid (x, y) \notin \Theta_{n} \lor (y, z) \notin \Theta_{n} } \\
  &= \set{ n \in \NN \mid (x, y) \notin \Theta_{n} } \cup \set{ n \in \NN \mid (y, z) \notin
    \Theta_{n} },
\end{align*}
we also have $d_{\Theta}(x, z) \leq \max\set{ d_{\Theta}(x, y), d_{\Theta}(y, z) }$ for
all $x, y, z \in X$.  Because of this $(X, d_{\Theta})$ is a pseudo-ultrametric space.

\begin{Proposition}
  \label{prop:entropic.dimension}
  Let $X$ be a non-empty set and let $\Theta = (\Theta_{n} \mid n \in \NN)$ be a
  descending sequence of equivalence relations on $X$ such that each $\Theta_{n}$ has
  finite index in $X$.  Then $(X, d_{\Theta})$ is precompact and
  \begin{equation*}
    \dim(X, d_{\Theta}) = \limsup_{n \to \infty} \frac{\log_{2}\abs{X / \Theta_{n}}}{n}.
  \end{equation*}
\end{Proposition}
\begin{Proof}
  We first observe that for all $x,y \in X$ and $n \in \NN$
  \begin{equation*}
    d_{\Theta}(x, y) < 2^{-n} \iff n < \inf\set{ m \in \NN \mid (x, y) \notin \Theta_{m} } \iff (x,
    y) \in \Theta_{n}.
  \end{equation*}
  Therefore, $X / \Theta_{n} = \set{ \Ball_{d_{\Theta}}(x, 2^{-n}) \mid x \in X }$.  Since $X /
  \Theta_{n}$ is finite, $(X, d_{\Theta})$ is precompact and
  \begin{equation*}
    \gamma_{(X, d_{\Theta})}(2^{-n}) = \abs{X / \Theta_{n}}.
  \end{equation*}
  Consequently,
  \begin{align*}
    \edim(X, d_{\Theta}) %
    &= \limsup_{r \to 0+} \frac{ \log_{2}(\gamma_{(X, d_{\Theta})}(r)) }{\log_{2}(1/r) }\\
    &= \limsup_{n \to \infty} \frac{ \log_{2}(\gamma_{(X, d_{\Theta})}(2^{-n})) }{ n }\\
    &= \limsup_{n \to \infty} \frac{ \log_{2} \abs{X / \Theta_{n} } }{ n }
  \end{align*}
  as required.
\end{Proof}

A straightforward application of this lemma is the following theorem.

\begin{Corollary}
  Let $\Sigma$ be an alphabet and let $L \subseteq \Sigma^{*}$.  Then with $\Theta :=
  (\Theta(\Sigma^{(n)}, L) \mid n \in \NN)$
  \begin{equation*}
    \edim(\Sigma^{*}, d_{\Theta}) = h(L).
  \end{equation*}
\end{Corollary}

In the case that the language $L$ is represented by a topological automaton we obtain the
following result.

\begin{Theorem}
  Let $\mathcal{A} = (X, \Sigma, \alpha, x_{0}, F)$ be a topological automaton. Let
  $\Lambda = (\Lambda_{n} \mid n \in \NN)$ where
  \begin{equation*}
    \Lambda_{n} \coloneqq \Lambda_{\Sigma^{(n)}} = \set{ (x, y) \in X \times X \mid \forall w \in
      \Sigma^{(n)} \colon \alpha(w,x) \in F \iff \alpha(w, y) \in F }
  \end{equation*}
  whenever $n \in \NN$ (cf.~\Cref{lem:entropy-of-automaton-as-entropy-of-cover}). Then
  $h(L(\mathcal{A})) \leq \edim(X, d_{\Lambda})$. Furthermore, if $\mathcal{A}$ is trim, then
  $h(L(\mathcal{A})) = \edim(X, d_{\Lambda})$.
\end{Theorem}

\begin{Proof}
  Let $L \coloneqq L(\mathcal{A})$ and $n \in \NN$. We observe that
  $\gamma_{L}(\Sigma^{(n)}) = \abs{\Sigma^{\ast}/\Theta(\Sigma^{(n)},L)} \leq
  X/\Lambda_{n}$
  by \Cref{lem:entropy-of-automaton-as-entropy-of-cover}~(2). Moreover, if $\mathcal{A}$
  is trim, then $\gamma_{L}(\Sigma^{(n)}) = X/\Lambda_{n}$ due to
  \Cref{lem:entropy-of-automaton-as-entropy-of-cover}~(3). Hence,
  \Cref{prop:entropic.dimension} yields the desired statements.
\end{Proof}

The pseudo-metric considered in the theorem above does not necessarily generate the topology of the
respective automaton. In fact, this happens to be true if and only if the automaton is minimal,
i.e., isomorphic to the minimal automaton of the accepted language. Furthermore, this case can be
characterized in terms of a separation property: a topological automaton is minimal if and only if
the induced pseudo-metric is a metric.

\begin{Proposition}
  Let $\mathcal{A} = (X, \Sigma, \alpha, x_{0}, F)$ be a topological automaton and
  $L \coloneqq L(\mathcal{A})$. Then the topology generated by $d_{\Lambda}$ is contained in the
  topology of $X$. Furthermore, the following statements are equivalent:
  \begin{enumerate}
  \item $\mathcal{A} \cong \mathcal{A}(L)$.
  \item $d_{\Lambda}$ is a metric.
  \item $d_{\Lambda}$ generates the topology of $X$.
  \end{enumerate}
\end{Proposition}

\begin{Proof}
  By \Cref{lem:entropy-of-automaton-as-entropy-of-cover}~(1), the subset
  $B_{d_{\Lambda}}(x,\epsilon) = [x]_{\Lambda_{-\lceil \log_{2}\epsilon \rceil}}$ is open in $X$ for
  all $x \in X$ and $\epsilon \in (0,\infty)$. Hence, the topology generated by $d_{\Lambda}$ is
  contained in the original topology of $X$. Now let us prove the claimed equivalences:

  (2)$\Longrightarrow$(3): Suppose that $d_{\Lambda}$ is a metric. Then the topology generated by
  $d_{\Lambda}$ is a Hausdorff topology. Since this topology is contained in the compact topology of
  $X$, both topologies coincide due to a basic result from set-theoretic topology
  (see~\cite[\S9.4, Corollary 3]{Bourbaki1}).

  (3)$\Longrightarrow$(1): Assume that $d_{\Lambda}$ generates the topology of $X$. This clearly
  implies $d_{\Lambda}$ to be a metric. Consider the unique surjective continuous homomorphism
  $\phi \colon \mathcal{A} \to \mathcal{A}(L)$. We are going to show that $\phi$ is injective. To
  this end, let $x,y \in X$ such that $\phi (x) = \phi (y)$. We argue that $d_{\Lambda}(x,y) =
  0$.
  Let $n \in \NN$. For every $w \in \Sigma^{(n)}$, we observe that
  \begin{align*}
    \alpha (x,w) \in F \
    &\Longleftrightarrow \ \phi (\alpha (x,w)) \in T_{L}
      \ \Longleftrightarrow \ \delta (\phi (x),w) \in T_{L} \\
    &\Longleftrightarrow \ \delta (\phi (y),w) \in T_{L}
      \ \Longleftrightarrow \ \phi (\alpha (y,w))
      \in T_{L} \ \Longleftrightarrow \ \alpha (y,w) \in F .
  \end{align*}
  Thus, $(x,y) \in \Lambda_{n}$. It follows that
  $(x,y) \in \bigcap_{n \in \NN} \Lambda_{n}$ and hence $d_{\Lambda}(x,y) = 0$. Since
  $d_{\Lambda}$ is a metric, we conclude that $x = y$. Accordingly, $\phi$ is a bijective
  continuous map between compact Hausdorff spaces and therefore a homeomorphism. This
  again is due to an elementary result from set-theoretic topology (see~\cite[\S9.4,
  Corollary 2]{Bourbaki1}).

  (1)$\Longrightarrow$(2): Suppose $\phi \colon \mathcal{A} \to \mathcal{A}(L)$ to be the
  necessarily unique isomorphism. Concerning any two points $x,y \in X$, we observe
  that
  \begin{align*}
    (x,y) \in \Lambda_{n}(\mathcal{A}) \
    &\Longleftrightarrow \ \forall w \in \Sigma^{(n)} \colon
      \, \alpha (x,w) \in F \Leftrightarrow \alpha (y,w) \in F  \\
    &\Longleftrightarrow \ \forall w \in \Sigma^{(n)} \colon
      \, \phi (\alpha (x,w)) \in T_{L} \Leftrightarrow \phi(\alpha (y,w)) \in T_{L}  \\
    &\Longleftrightarrow \ \forall w \in \Sigma^{(n)} \colon
      \, \delta (\phi(x),w) \in T_{L} \Leftrightarrow \delta (\phi(y),w) \in T_{L}  \\
    &\Longleftrightarrow \ (\phi(x),\phi(y)) \in
      \Lambda_{n}(\mathcal{A}(L))
  \end{align*}
  for every $n \in \NN$. Hence,
  $d_{\Lambda (\mathcal{A})}(x,y) = d_{\Lambda (\mathcal{A}(L))}(\phi(x),\phi(y))$ for all
  $x,y \in X$. Accordingly, it suffices to show that $d_{\Lambda (\mathcal{A}(L))}$ is a metric. To
  this end, let $f,g \in \overline{\chi_{L}(\Sigma^{\ast})}$ such that
  $d_{\Lambda (\mathcal{A}(L))} (f,g) = 0$. We argue that $f = g$. Let $w \in \Sigma^{\ast}$.  Then
  there exists $n \in \NN$ where $w \in \Sigma^{(n)}$. Since
  $d_{\Lambda (\mathcal{A}(L))} (f,g) = 0$, we conclude that $(f,g) \in \Lambda_{n}(\mathcal{A}(L))$
  and thus
  \begin{align*}
    f(w) = 1 \
    &\Longleftrightarrow \ \delta (f,w)(\epsilon) = 1 \ \Longleftrightarrow \ \delta (f,w) \in T_{L} \\
    &\Longleftrightarrow \ \delta (g,w) \in T_{L} \ \Longleftrightarrow \ \delta
      (g,w)(\epsilon) = 1 \ \Longleftrightarrow \ g(w) = 1 .
  \end{align*}
  Therefore, $f(w) = g(w)$. It follows that $f=g$. This shows that $d_{\Lambda (\mathcal{A}(L))}$ is
  a metric and hence completes the proof.
\end{Proof}

\section{Conclusions}
\label{sec:conclusions}

In this paper we have introduced the notion of topological entropy of formal
languages as the topological entropy of the minimal topological automaton
accepting it.  We have shown that this notion is equal to the Myhill-Nerode
complexity of the language, and can also be characterized in terms of the
entropic dimension of suitable pseudo-ultrametric spaces.  Using these
characterizations, we were able to compute the topological entropy of certain
types of languages.

The main motivation of this work was the goal to uniformly assess the complexity
of formal languages independent of a particular collection of computation
models.  We believe that the examples we have provided in this work show that
the notion of topological entropy of formal languages is a suitable candidate
for such a complexity measure.  In particular, we have shown that some languages
intuitively considered to be simple all have zero entropy: regular languages,
Dyck languages with one sort of parentheses, our \enquote{commutative} version
of Dyck languages with arbitrary sorts of parentheses, and the language $\set{
  a^{n}b^{n}c^{n} \mid n \in \mathbb N }$.  Indeed, all of these languages are
accepted by simple models of computation, e.g., one-way finite automata with a
fixed number of counters.

On the other hand, we have presented examples of languages that have non-zero
entropy that can hardly be considered as simple, namely Dyck languages with more
than one sort of parentheses as well as the palindrome languages.  Indeed,
palindromes cannot be accepted by deterministic pushdown automata, and Dyck
languages with more than one sort of parentheses give rise to the hardest
context-free languages~\cite{book/HoFL/ContextFreeLanguages}.

A natural next step in investigating the notion of topological entropy is to
provide more examples that test the suitability of this notion as a measure of
complexity of formal languages.  For example, we have already shown that all
languages accepted by finite automata have zero entropy.  A natural question is
now to ask for which classes of computation models the topological entropy of
the accepted languages is also zero.  We suspect that this is the case for
one-way finite automata equipped with a fixed number of counters and an
acceptance condition that does only require to check local conditions, including
the current values of the counters.

Conversely, one could ask what properties languages with non-zero entropy
possess.  What form of non-locality in a suitable machine model is necessary to
accept such languages, given that they are decidable?  And what properties do
languages have if their topological entropy is infinite?  Are there context-free
languages with infinite entropy?

\bigskip

\textit{Acknowledgments}: We would like to express our sincere gratitude towards
the anonymous reviewer for numerous valuable suggestions that significantly
improved the presentation of the paper.

\appendix

\printbibliography

\end{document}